%%%%%%%%%%%%%%%%%%%% This is a Plain-Tex file.
%**start of header
\magnification=1200

\overfullrule=0pt

\abovedisplayskip=8pt plus 3pt minus 1pt
\belowdisplayskip=8pt plus 3pt minus 1pt
\abovedisplayshortskip=2pt plus 3pt
\belowdisplayshortskip=6pt plus 3pt

% CALL 12pt FONTS
\font\twelverm=cmr12
\font\twelvebf=cmbx12
\font\twelveit=cmti12
\def\bigfont{		% DEFINE 12pt FONTS
	\let\rm=\twelverm
	\let\bf=\twelvebf
	\let\it=\twelveit
	\rm}

\input amssym.def
\input amssym
% for Letter size
\vsize=8truein
\hsize=6truein
% for A4 siz
%\vsize=24true cm
%\hsize=16true cm
\hoffset=18truept

\catcode`\@=11
\newskip\ttglue
\font\ninerm=cmr9

\font\sixrm=cmr6

\font\ninei=cmmi9
\font\eighti=cmmi8
\font\sixi=cmmi6
\skewchar\ninei='177 \skewchar\eighti='177 \skewchar\sixi='177

\font\ninesy=cmsy9
\font\eightsy=cmsy8
\font\sixsy=cmsy6
\skewchar\ninesy='60 \skewchar\eightsy='60 \skewchar\sixsy='60

\font\ninebf=cmbx9

\font\sixbf=cmbx6

\font\ninett=cmtt9
\font\eighttt=cmtt8

\hyphenchar\tentt=-1 % inhibit hyphenation in typewriter type
\hyphenchar\ninett=-1
\hyphenchar\eighttt=-1

\font\ninesl=cmsl9

\font\nineit=cmti9

\def\ninepoint{\def\rm{\fam0\ninerm}%
  \textfont0=\ninerm \scriptfont0=\sixrm \scriptscriptfont0=\fiverm
  \textfont1=\ninei \scriptfont1=\sixi \scriptscriptfont1=\fivei
  \textfont2=\ninesy \scriptfont2=\sixsy \scriptscriptfont2=\fivesy
  \textfont3=\tenex \scriptfont3=\tenex \scriptscriptfont3=\tenex
  \def\it{\fam\itfam\nineit}%
  \textfont\itfam=\nineit
  \def\sl{\fam\slfam\ninesl}%
  \textfont\slfam=\ninesl
  \def\bf{\fam\bffam\ninebf}%
  \textfont\bffam=\ninebf \scriptfont\bffam=\sixbf
   \scriptscriptfont\bffam=\fivebf
  \def\tt{\fam\ttfam\ninett}%
  \textfont\ttfam=\ninett
  \tt \ttglue=.5em plus.25em minus.15em
  \normalbaselineskip=11pt
  \def\MF{{\manual hijk}\-{\manual lmnj}}%
  \let\sc=\sevenrm
  \let\big=\ninebig
  \setbox\strutbox=\hbox{\vrule height8pt depth3pt width\z@}%
  \normalbaselines\rm}

\font\smaller=cmr8
\font\ninecsc=cmcsc9
\font\tencsc=cmcsc10

\def\beginsection#1{\bigskip\medskip
\centerline{\tencsc #1}\bigskip}

\catcode`\@=12

\headline={\ifnum\pageno>1\smaller\ifodd\pageno\hfill 
NEVANLINNA THEORY AND RATIONAL POINTS
\hfill\the\pageno \else
\the\pageno\hfill\uppercase{Junjiro Noguchi}\hfill\fi\else\hss\fi}

\footline={\hss}

\centerline{\bf ON HOLOMORPHIC CURVES IN SEMI-ABELIAN VARIETIES}
\footnote{}{Research at MSRI supported in part by NSF grant \#DMS 9022140.}
\bigskip\medskip
\centerline{\bf Junjiro Noguchi}

\midinsert
\narrower\narrower
\noindent {\ninecsc Abstract.}
\baselineskip=10pt
{\ninepoint 
The algebraic degeneracy of holomorphic curves in a semi-Abelian
variety omitting a divisor is proved
(Lang's conjecture generalized to semi-Abelian varieties)
by making use of the {\it jet-projection method} and the
logarithmic Wronskian jet differential after Siu-Yeung.
We also prove a structure theorem for the locus which contains
all possible image of non-constant entire holomorphic curves
in a semi-Abelian variety omitting a divisor.

}
\endinsert

\beginsection{\indent Introduction}

\indent
The purpose of this paper is to prove this result:

\medskip
\noindent
{\bf Main Theorem.}
{\it
Let $D$ be a non-zero algebraic effective reduced divisor of a semi-Abelian
variety $A$ over the complex number field $\bf C$.
Let $f:{\bf C} \rightarrow A \setminus D$ be an arbitrary holomorphic mapping.
}
\parindent=30pt

\item{{\rm (i)}}
{\it
The Zariski closure $X_0(f)$ of the image of $f$ in $A$ is a translate of
a proper semi-Abelian subvariety of $A$, and $X_0(f)\cap D=\emptyset$.
}
\item{\rm (ii)}
{\it
In special, if $D$ has a non-empty intersection with
any translate of any positive dimensional semi-Abelian subvariety of $A$,
then $f:{\bf C} \rightarrow A \setminus D$ is constant.
}
\parindent=10pt

\medskip
Moreover, we prove a {\it structure theorem} for the locus (sometimes, called
the exceptional set) of $A\setminus D$ which
contains the images of all possible non-constant entire holomorphic curves
in $A\setminus D$ (see Theorem (3.1) and Remark after it).

In the case where $A$ is a semi-Abelian variety and $D$ has two components
which are homologous to each other, the algebraic degeneracy of
a holomorphic curve $f:{\bf C} \rightarrow A \setminus D$ was proved by [N81].
In the case where $A$ is an Abelian variety, Siu-Yeung [SY96] proved the
Main Theorem (Lang's conjecture).
Note that a generalization of the Main Theorem to the case of
semi-Abelian varieties is already claimed in the introduction of [SY96],
and that the lemmas we are going to prove are similar to those in [SY96].
The main difference is that instead of Siu-Yeung's elaborate
Wronskian arguments ([SY96], Lemma (1.1) and Proof of Lemma (2.1)),
we use a simpler and more direct ``{\it jet-projection method\/}'',
the same idea as in [NO${84}\over{90}$], Chap.\ VI; there, a self-contained and
detailed proof of Bloch's conjecture is described.
This proof, especially the proof of [NO${84}\over{90}$], Lemma (6.3.10), is an original
one by the present author for the remaining part after [B26] and [O77]
to finish the proof of Bloch's conjecture,
and works as well for semi-Abelian varieties [N81].
M. Green gave a talk on a proof of [NO${84}\over{90}$], Lemma (6.3.10) based on
Gauss' maps at Taniguchi Symposium, Katata/Kyoto, 1978;
he did not publish it, but did [GG80] with P. Griffiths by making use of
another idea based on Riemann-Roch and curvature methods
(cf.\ also [K80] and [M96]);
afterward, a gap in the proof of [GG80] was found;
see  [D96] and [DL96], where Dethloff and Lu also use the present
{jet-projection method}.

It is worth noting that Ochiai [O77] proved the algebraic degeneracy of
an entire holomorphic curve in an Abelian variety omitting two divisors,
mutually linearly equivalent, and recognizing that this is the first link
between Bloch's and Lang's conjectures.

In the course of the proof of the Main Theorem,
we also need a generalization of theta functions
for semi-Abelian varieties (see Lemma (2.1) and Remark after it).

In comparison with arithmetic on semi-Abelian varieties,
the counter part in number theory was proved by Vojta [V95].

{\it Acknowledgement.}
The present paper was written during the author's stay at
Mathematical Sciences Research Institute, University of California, Berkeley
for the special year program in Several Complex Variables, 1995/96.
The author is grateful to the institute for the very active
and stimulating circumstance, and to participating members,
especially, Professors Y.-T.\ Siu, B. Shiffman,
S. Lu, G. Dethloff, M. Ru, P.-M. Wong, P. Vojta, A. Huckleberry,
F. Catanese, H. Tsuji, and S. Kobayashi.

\beginsection{\indent\S1. Jet space and translation invariance}

\indent
Let $A$ be a semi-Abelian variety of dimension $n$ over $\bf C$;
that is, $A$ is a complex algebraic group of dimension $n$ which
is an extension of an Abelian variety $A_0$ by an algebraic
torus $({\bf C}^*)^t$ ($t \in {\bf Z}, \geqq 0$),
$$
0 \rightarrow ({\bf C}^*)^t \rightarrow A \rightarrow A_0 \rightarrow 0.
\leqno(1.1)
$$
(Cf.\ [I76].)
Let $f:{\bf C} \rightarrow A$ be a holomorphic curve.
Let $J_k(A)$ denote the $k$-th jet space over $A$ with $k \geqq 0$
(where if $k=0$, $J_0(A)=A$), and let
$d^k f:{\bf C} \rightarrow J_k(A)$ denote the $k$-th jet lifting of $f$.
(Cf.\ [NO${84}\over{90}$], Chap.\ VI.)
Let $X_k(f)$ be the Zariski closure of $d^k f({\bf C})$ in $J_k(A)$.

Note that the group law of $A$ canonically induces an additive
action of $A$ on $J_k(A)$.

{\bf Lemma (1.2).}
{\it If $\dim X_k(f) >0$, then there exists a one-parameter subgroup
of $A$ which leaves $X_k(f)$ invariant.}

{\it Proof.}
Since $A$ is a semi-Abelian variety, we have a decomposition
$$
J_k(A)=A \times {\bf C}^{nk},
\leqno(1.3)
$$
which is fixed from now on.

We take the $l$-th jet space $J_l(X_k(f))$ (resp.\ $J_l(J_k(A))$)
of $X_k(f)$ (resp.\ $J_k(A)$), where $l=1,2,3, \ldots$.
Then $J_l(X_k(f)) \subset J_l(J_k(A))$, and we obtain
the decomposition naturally induced from (1.3)
$$
J_l(J_k(A))=J_l(A \times {\bf C}^{nk})=
A \times {\bf C}^{nk} \times {\bf C}^{(n+nk)l},
\leqno(1.4)
$$
the natural projection to the second and third factors, and its restriction:
$$
\eqalign{
\psi_{k,l} &:J_l(J_k(A)) \rightarrow {\bf C}^{nk}\times{\bf C}^{(n+nk)l},\cr
\Psi_{k,l}=\psi_{k,l}|J_l(X_k(f))&:J_l(X_k(f)) \rightarrow
{\bf C}^{nk}\times {\bf C}^{(n+nk)l}.}
\leqno(1.5)
$$
Let $d^ld^kf:{\bf C} \rightarrow J_l(X_k(f))$ be the $l$-the
jet lifting of $d^kf$.
Take a point $z_0 \in {\bf C}$ with $d^kf(z_0) \in X_k(f)_{\hbox{reg}}$,
where $X_k(f)_{\hbox{reg}}$ denotes the set of regular points of $X_k(f)$.
Set $y_l=d^ld^kf(z_0) \in J_l(X_k(f))$.
Now, look at the kernel
$$
\eqalign{
\hbox{Ker }d\Psi_{k,l}(y_l)
&
\subset
{\bf T}_{y_l} (J_l(X_k(f))) \subset
{\bf T}_{y_l} (J_l(J_k(A)))
\cr
&
={\bf T}_{f(z_0)}(A) \oplus
{\bf T}_{\psi_{k,l}(y_l)}({\bf C}^{nk} \times {\bf C}^{(n+nk)l})
}
\leqno(1.6)
$$
(see (1.4))
of the differential $d\Psi_{k,l}$ of $\Psi_{k,l}$ at the point
$y_l$,
where ${\bf T}(*)$ stands for the holomorphic tangent space.
By the definition of $\Psi_{k,l}$,
$$
\hbox{Ker }d\Psi_{k,l}(y_l) \subset {\bf T}_{f(z_0)}(A)
\oplus O' \cong{\bf T}_{f(z_0)}(A),
\leqno(1.7)
$$
where $O'$ denotes the zero vector of
${\bf T}_{\psi_{k,l}(y_l)}({\bf C}^{nk} \times {\bf C}^{(n+nk)l})$
in (1.6).
By making use of the isomorphism in (1.7),
$\hbox{Ker }d\Psi_{k,l'}(y_{l'}) \subset
\hbox{Ker }d\Psi_{k,l}(y_l) \subset {\bf T}_{f(z_0)}(A)$
for $l' \geqq l$.
Therefore, if
$\cap_{l \geqq 1} \hbox{Ker }d\Psi_{k,l}(y_l) =\{O\}$,
$\hbox{Ker }d\Psi_{k,l}(y_l)=\{O\}$ for some $l$.
Then, by making use of Lemma on logarithmic derivative,
we obtain an estimate of the order function $T_f(r)$ of $f$
by those of $df, d^2f, \ldots, d^{k+l}f$ which are
small functions with respect to $T_f(r)$;
this is a contradiction (cf.\ [NO${84}\over{90}$], Chap.\ VI
for this argument).
Thus, $\cap_{l \geqq 1} \hbox{Ker }d\Psi_{k,l}(y_l) \not =\{O\}$.
Let $v \in \cap_{l \geqq 1} \hbox{Ker }d\Psi_{k,l}\setminus \{O\}
\subset {\bf T}_{f(z_0)}(A)\setminus\{O\}$.
We may regard $v$ as a global holomorphic vector field on $J_k(A)$ through the
trivialization (1.3).
Then $v$ is tangent to $X_k(f)$ at all points of the image of a
neighborhood of $z_0$ by $d^kf$, and so it is tangent to $X_k(f)$
at all points of $X_k(f)$.
Hence $X_k(f)$ is invariant by the action of the one-parameter
subgroup generated by $v$.   {\it Q.E.D.}

We define the stabilizer ${\rm St}(X_k(f))$ to be the maximal semi-Abelian
subvariety of $A$ which leaves $X_k(f)$ invariant by translations.
After taking the quotient by ${\rm St}(X_k(f))$, we apply Lemma (1.2)
to conclude

{\bf Proposition (1.8).}
{\it 
Let the notation be as above.}
\parindent=30pt

\item{\rm (i)}
{\it
$X_0(f)$ is a translation of a semi-Abelian subvariety $B$ of $A$.
}
\item{\rm (ii)}
{\it
$X_k(f)$ is invariant by the action of $B$ on $J_k(A)$ for all $k \geqq 0$.
}
\parindent=10pt

\beginsection
{\indent\S2. Proof of Main Theorem}
 
By (1.1) we consider $A$ as a $({\bf C^*})^t$-principal fiber bundle
with projection $\rho:A \rightarrow A_0$.
Using the natural compactification $({\bf C^*})^t \subset ({\bf P}^1_{\bf C})^t$,
we have a compactification
$\overline A$ of $A$ with projection $\rho:\overline A \rightarrow A_0$,
which is projective algebraic (cf.\ [N81] and [V95]).
Since $D$ be an algebraic effective reduced
divisor on $A$, there is an
effective reduced divisor $\overline D$ on $\overline A$
with $\overline D \cap A= D$.

Let $\pi:{\bf C}^{n} \rightarrow A$ be the universal covering of $A$ with
$\pi_1(A)=\Gamma$ ([I76]), which may be called an incomplete lattice or
a semi-lattice.

{\bf Lemma (2.1).}
{\it
Let the notation be as above.
Then there is an entire function $\theta(x)$ on ${\bf C}^n$
such that $\pi^*D=(\theta)$ and for $\gamma \in \Gamma$
$$
\theta(x+\gamma)=
e^{L_\gamma(x)} \theta(x),\qquad x \in {\bf C}^{n},
$$
where $L_\gamma(x)$ is affine linear in $x$.
}

{\it Remark.}
In a glance this Lemma (2.1) seems to be already known and classical,
but so far, asking several specialists, we could not find
a reference stating the above lemma.
So we give a proof, which is actually quite elementary.

{\it Proof.}
We may assume that $\overline D$ is irreducible, and set $n_0=\dim A_0$.
If $\rho (\overline D) \not=A_0$, then $\dim \rho (\overline D)=n_0-1$,
and $\rho^{-1}(\rho (\overline D))\cap A=D$.
Then this is the case of the classical theta function on an Abelian
variety (cf.\ [W58]).

Assume that $\rho (\overline D) =A_0$.
Let $\{U_{\lambda}\}$ be an open covering of $A_0$ such that
$U_\lambda$ is a ball in ${\bf C}^{n_0}$ ($n_0=\dim A_0$),
and the restriction $A|U_\lambda$ is isomorphic to
$U_\lambda \times ({\bf C^*})^t$.
We fix a multiplicative coordinate system
$(z_1, \ldots, z_t)$ of $({\bf C^*})^t$.
Then $D\cap (A|U_\lambda)$ is given by a zero set of
a finite Laurent series in $(z_1, \ldots, z_t)$
$$
\theta_\lambda=\sum_{\rm finite}
a_{\lambda l_1 \cdots l_t}(y) z_1^{l_1} \cdots z_t^{l_t},
$$
where the coefficients $a_{\lambda l_1 \cdots l_t}(y)$ are holomorphic
functions in $U_\lambda$.
Then, on $(A|U_\lambda) \cap (A|U_\mu)$ we have
$$
\theta_\lambda \cdot \theta_\mu^{-1}=a_{\lambda\mu}(y)
z_1^{l_{\lambda\mu 1}} \cdots z_t^{l_{\lambda\mu t}},
$$
where the decomposition of the right side is unique.
Note that $\{a_{\lambda\mu}(y)\}$ forms a 1-cycle of non-vanishing
holomorphic functions on $A_0$, and hence defines a holomorphic line bundle
over $A_0$.
Making use the classical theta function on $A_0$, we obtain
the required entire function $\theta$.
{\it Q.E.D.}

We take $\theta(x)$ in Lemma (2.1), which may be called the
theta function associated with $D$.
Choose and fix a linear coordinate system $(x_1, \ldots, x_n)$ of
${\bf C}^{n}$,
so that the lifting $\tilde f: {\bf C}\rightarrow {\bf C}^{n}$
of $f$ is expressed as
$$
\tilde f(z)=(f_1(z), \ldots, f_{n'}(z), 0,\ldots,0)
$$
with entire functions, $f_1(z), \ldots, f_{n'}(z)$, which are linearly
independent over ${\bf C}$.
We define an algebraic
jet subbundle $J_k(A)'$ of $J_k(A)$ ($k=0,1,\ldots$) by equations
$$
d^i x_j=0,\qquad 1 \leqq i \leqq k,\quad n'+1 \leqq j \leqq n
$$
(cf.\ the notion of directed manifolds in [D96]).
Then, $X_{n'+1}(f) \subset J_{n'+1}(A)'$, and
after Siu-Yeung [SY96] we define the logarithmic jet differential on
$J_{n'+1}(A)'$
$$
\Theta=\left| \matrix{
d\log \theta & dx_1 & \cdots & dx_{n'} \cr
\vdots & \vdots & \vdots & \vdots \cr
d^{n'+1}\log \theta & d^{n'+1}x_1 & \cdots & d^{n'+1}x_{n'} \cr
}\right| ,
\leqno{(2.2)}
$$
which is well-defined by Lemma (2.1).
Denote by the same $\Theta$ the restriction of $\Theta$ over
$X_{n'+1}(f)$:
$$
\Theta : X_{n'+1}(f) \rightarrow {\bf C},
$$
which is a rational function with logarithmic poles on fibers over $D$.
Then, taking the derivatives of $\Theta$, we have
$$
d^l\Theta : J_l(X_{n'+1}(f)) \rightarrow {\bf C},
\qquad l=0,1,2,\ldots.
$$
Set
$$
\eqalign{
\Lambda_{n'+1,l} &=\Psi_{n'+1,l} \times
(\Theta\circ p_{0, l}, d\Theta\circ p_{1,l},\ldots,
d^l\Theta \circ p_{l,l})\cr
&: J_l(X_{n'+1}(f))\rightarrow
{\bf C}^{n(n'+1)}\times {\bf C}^{(n+n(n'+1))l} \times {\bf C}^{l+1},}
$$
where $p_{j, l}:J_l(X_{n'+1}(f)) \rightarrow J_j(X_{n'+1}(f))$,
$0 \leqq j \leqq l$, are the canonical projections.
Then we apply the same argument as in the proof of Lemma (1.2)
for $\Lambda_{n'+1,l}$ in place of $\Psi_{k,l}$ to
deduce

{\bf Lemma (2.3).}
{\it If $\dim X_{n'+1}(f) >0$, then there is a non-zero holomorphic
vector field $v'$ on $A$ such that $X_k(f)$, $k \geqq 0$,
are invariant by the translations generated by $v'$,
and $v'\Theta\equiv 0$ on $X_{n'+1}(f)$.
}

Now, assume that $f:{\bf C} \rightarrow A \setminus D$ is not
constant.
Then, $\dim X_{n'+1}(f) >0$ for some $n'\geq 1$, and we may use Lemma (2.3).
Since $v'\Theta (d^{n'+1}f(z))\equiv 0$ in $z$,
there are complex numbers $c_1$, $\ldots$, $c_{n'}$ such that
$$
d v'\log \theta (f(z))+ c_1 df_1(z)+ \cdots + c_{n'} df_{n'}(z)\equiv 0.
$$
Therefore
$$
d v'\log \theta + c_1 dx_1+ \cdots + c_{n'} dx_{n'}\equiv 0
$$
on $X_1(f)$.
By Lemma (2.3), $X_1(f)$ is invariant by the translations $\hbox{Exp}(t v')$,
$t \in {\bf C}$, and so $d v' v'\log \theta \equiv 0$ on $X_1(f)$.
Thus, $d v' v'\log \theta(f(z)) \equiv 0$, and then
$$
v' v'\log \theta \equiv \hbox{constant on } X_0(f).
\leqno(2.4)
$$
Note that by Proposition (1.8) and the assumptions,
$X_0(f)$ is a translation of
a proper semi-Abelian subvariety $B$ of positive dimension, and that
$v'$ is tangent to $B$ at all points of $B$.
It follows from (2.4) that $X_0(f)\setminus D$ and $X_0(f)\cap D$ are
invariant by the translations generated by $v'$.

Suppose that $X_0(f)\cap D\not=\emptyset$, and set $D'=X_0(f)\cap D$.
By a translation we may assume that $B=X_0(f)$.
Then $D'$ is a non-zero divisor of $B$.
Let ${\rm St}(D') \subset B$ be the stabilizer of $D'$ in $B$.
Let $g: {\bf C} \rightarrow (X_0/{\rm St}(D'))\setminus (D'/{\rm St}(D'))$
be the composition
of $f$ and the quotient mapping by the action of ${\rm St}(D')$.
Then $g$ is not algebraically degenerate.
We apply again the above proved for $g$ to deduce the algebraic degeneracy
of $g$. This is a contradiction.
So we have proved (i) and (ii).

This completes the proof of the Main Theorem.

\beginsection
{\indent\S3. Translates of semi-Abelian subvarieties in $A \setminus D$}

\indent
Let $A$ and $D$ be as in the Main Theorem.
Let $X$ be an algebraic subset of $A$, and
let $E(X,D)$ denote the set of all translates of semi-Abelian subvarieties
of $A$ which are contained in $X$ and have no intersection with $D$.

{\bf Theorem (3.1).}
{\it
Let the notation be as above.
Then $E(X,D)$ is an algebraic subset such that it decomposes to
irreducible components $E_i$ with $\dim {\rm St}(E_i) > 0$.
}

{\it Remark.}
In the case where $A$ is an Abelian variety and $D=\emptyset$,
this was proved by Kawamata [K80], and when $A$ is a semi-Abelian
variety and $D=\emptyset$, it was proved by [N81].
Vojta [V95] generalized it to the case of $A$ and $D$ as in
the above Theorem (3.1).
By making use of the results of \S\S1 and 2, we are going to give another
proof, which is totally different to Vojta's [V95].
Abramovich [A94] gave also a proof for the case of $D=\emptyset$,
which works over fields of arbitrary characteristic $\geqq 0$.
The way of Abramovich [A94] referring to the result in [N81]
would be misleading.

{\it Proof.}
Let $\pi:{\bf C}^{n} \rightarrow A$ be the universal covering of $A$ with
semi-lattice $\Gamma$.
Let ${\bf P}({\bf C}^n)$ be the projective space of 1-dimensional linear
subspaces $[v]$ $(={\bf C} v)$ with  $v \in {\bf C}^n\setminus\{O\}$.
Set
$$
\eqalign{
      X'&=X\setminus D\cr
{\cal E}&=\{(x,[v]) \in X' \times {\bf P}({\bf C}^n);
x+{\bf C} v \subset X'\}, \cr
\mu &: (x,[v]) \in {\cal E} \rightarrow  x \in X'.}
$$
Since $D$ is a Cartier divisor, $\cal E$ is closed in
$X' \times {\bf P}({\bf C}^n)$, so that $\mu$ is proper.
We use $\Psi_{0,l}$ in (1.5) with taking $X_0(f)=X$, and $\theta$ in (2.4).
It follows from the proofs of Lemmas (1.2) and (2.3) that
$\cal E$ is the set of points $(x,[v]) \in X' \times {\bf P}({\bf C}^n)$
defined by algebraic equations
$$
v \in \cap_{l\geqq 1} \hbox{Ker }d\Psi_{0,l}(x),\quad
\hbox{and} \quad
v^k\log \theta (x)=0,\quad k \geqq 3,
$$
where $v$ is identified with a holomorphic vector field on $A$ and
$v^k$ stands for the $k$-th derivation in the direction $v$.
Thus $\cal E$ is algebraic, and so is
$E(X,D)=\mu({\cal E})$ in $X'$.

By making use of the countability of semi-Abelian subgroups of $A$ and
the Baire's category theorem we infer that $\dim {\rm St}(E_i) > 0$
for every irreducible component $E_i$ of $E(X,D)$
(cf.\ the proof of [N81], Lemma (4.1)).
{\it Q.E.D.}

{\it Example.} 1) Let $A$ be an Abelian variety and let $D$ be a divisor.
Then $E(A,D) \subsetneqq A$ if and only if $D$ is ample;
moreover, if $D$ is ample, $E(A,D)=\emptyset$ (cf., e.g., [W58]).

2) Let $A=({\bf C}^*)^t$ with coordinates $(z_1,\ldots,z_t)$ and let $D$ be
defined by
$$
z_1 + \cdots + z_t=1.
$$
Let $f:{\bf C} \rightarrow A\setminus D$ be a holomorphic curve.
It follows from Borel's theorem that after a change of order of coordinates
$$
\matrix{
f_1 \sim \cdots \sim f_{j_1}\sim 1 & \hbox{and} & f_1+\cdots +f_{j_1}=1 \cr
f_{j_1+1} \sim \cdots \sim f_{j_2} & \hbox{and} & f_{j_1+1}
+\cdots +f_{j_2}=0 \cr
\cdots & \cdots & \cdots \cr
f_{j_k+1} \sim \cdots \sim f_{n} & \hbox{and} & f_{j_k+1}
+\cdots +f_{n}=0,}
$$
where $f_1 \sim f_2$ means the constancy of $f_1/f_2$.
Thus $E(A,D)$ consists of this kind of translates of subgroups.
%; e.g., for $t=2$,

%\vskip 40true mm
%\centerline{Figure}

See [N81] and [V95] for more examples.

\beginsection{References}

\ninepoint

\item{[A94]}
D. Abramovich, Subvarieties of semiabelian varieties, Compositio Math.
{\bf 90} (1994), 37-52.
\item{[B26]}
A. Bloch,
Sur les syst\`emes de fonctions uniformes satisfaisant \`a l'\'equation d'une vari\'et\'e alg\'ebrique dont l'irr\'egularit\'e d\'epasse la dimension,
J. Math.\ Pures Appl.\
{\bf 5} (1926), 9-66.
\item{[D96]}
J.-P.\ Demailly,
Algebraic criteria for Kobayashi hyperbolic projective varieties and jet
differentials, preprint, 1996.
\item{[DL96]}
G. Dethloff and S. Lu, a seminar talk at MSRI, UC, Berkeley, 1996.
\item{[GG80]}
M. Green and P. Griffiths,
Two applications of algebraic geometry to entire holomorphic mappings,
The Chern Symposium 1979, pp.\ 41-74, Springer-Verlag,
New York-Heidelberg-Berlin,  1980.
\item{[I76]}
S. Iitaka,
Logarithmic forms of algebraic varieties,
J. Fac.\ Sci.\ Univ.\ Tokyo, Sect.\ IA {\bf 23} (1976), 525-544.
\item{[K80]}
Y. Kawamata,
On Bloch's conjecture, Invent.\ Math.\
{\bf 57} (1980),  97-100.
\item{[Ko96]}
S. Kobayashi,
Hyperbolic Complex Spaces, manuscript, 1996.
\item{[M96]}
M. McQuillan,
A new proof of the Bloch conjecture,
J. Algebraic Geometry
{\bf 5} (1996), 107-117.
\item{[N81]}
J.\ Noguchi,
Lemma on logarithmic derivatives and holomorphic curves in algebraic varieties,
Nagoya Math.\ J.\ {\bf 83} (1981), 213-233.
\item{[NO${84}\over{90}$]}
J. Noguchi and T. Ochiai,
Geometric Function Theory in Several Complex Variables,
Japanese edition, Iwanami, Tokyo, 1984;
English Translation, Transl.\ Math.\ Mono.\ {\bf 80},
Amer.\ Math.\ Soc., Providence, Rhode Island,
1990.
\item{[O77]}
T. Ochiai,
On holomorphic curves in algebraic varieties with ample irregularity,
Invent.\ Math.\ {\bf 43} (1977), 83-96.
\item{[SY96]}
Y.-T.\ Siu and S.-K.\ Yeung,
A generalized Bloch's theorem and the hyperbolicity of the complement
of an ample divisor in an Abelian variety,
preprint, 1996.
\item{[V95]}
P. Vojta, Integral points on subvarieties of semiabelian varieties, II,
preprint, 1995.
\item{[W58]}
A. Weil,
Introduction \`a l'\'Etude des Vari\'et\'es k\"ahleriennes,
Hermann, Paris, 1958.
\bigskip

\bigskip\rm
\noindent{\ninecsc
Department of Mathematics, Tokyo Institute of Technology, Ohokayama, Meguro,
Tokyo 152, Japan ({\tt noguchi@math.titech.ac.jp}).}

\medskip{\it \noindent Current address:}
Mathematical Sciences Research Institute, University of California, Berkeley,
1000 Centennial Drive, Berkeley, California 94720, USA
({\tt noguchi@\allowbreak msri.org}).

\bye